\documentclass[11pt,a4paper]{article}%
\usepackage{amscd}
\usepackage{amsmath}
\usepackage{graphicx}
\usepackage{amsfonts}
\usepackage{amssymb}%
\newtheorem{theorem}{Theorem}[section]
\newtheorem{lemma}[theorem]{Lemma}

\newtheorem{proposition}[theorem]{Proposition}

\newtheorem{remark}[theorem]{Remark}

\numberwithin{equation}{section}

\begin{document}

\title{Global well-posedness and symmetries for dissipative active scalar equations
with positive-order couplings}
\author{\textbf{Lucas C. F. Ferreira}\\{\small Universidade Estadual de Campinas, IMECC- Departamento de
Matem\'atica,} \\{\small Rua S\'ergio Buarque de Holanda, 651, CEP 13083-859, Campinas-SP,
Brazil.} \\{\small \texttt{email:\ lcff@ime.unicamp.br}} \vspace{1cm}\\\textbf{Lidiane S. M. Lima}\\{\small Universidade Estadual de Campinas, IMECC- Departamento de
Matem\'atica,} \\{\small Rua S\'{e}rgio Buarque de Holanda, 651, CEP 13083-859, Campinas-SP,
Brazil.} \\{\small \texttt{email:\ lidynet@hotmail.com}}}
\date{}
\maketitle

\begin{abstract}
We consider a family of dissipative active scalar equations outside the
$L^{2}$-space. This was introduced in [D. Chae, P. Constantin, J. Wu, \emph{to
appear in} IUMJ (2014)] and its velocity fields are coupled with the active
scalar via a class of multiplier operators which morally behave as derivatives
of positive order. We prove global well-posedness and time-decay of solutions,
without smallness assumptions, for initial data belonging to the critical
Lebesgue space $L^{\frac{n}{2\gamma-\beta}}(\mathbb{R}^{n})$ which is a class
larger than that of the above reference. Symmetry properties of solutions are
investigated depending on the symmetry of initial data and coupling operators.

\bigskip\noindent\textbf{AMS MSC:} 35Q35, 76D03, 35A01, 35B06, 35B40, 35R11, 86A10

\medskip\noindent\textbf{Keywords:}{ Active scalar equations, Global
well-posedness, Decay of solutions, Symmetry, {\nobreak{Critical}} spaces}

\end{abstract}


\pagestyle{myheadings} \markright{Dissipative active scalar equations}

\section{Introduction}

We are concerned with the initial value problem (IVP) for a family of
dissipative active scalar equation, which reads as
\begin{equation}%
\begin{cases}
\frac{\partial\theta}{\partial t}+\kappa\left(  -\Delta\right)  ^{\gamma
}\theta+u\cdot\nabla_{x}\theta=0, & \qquad x\in\mathbb{R}^{n}%
\ ,\ t>0,\ \\[3mm]%
\theta(x,0)=\theta_{0}(x), & \qquad x\in\mathbb{R}^{n},\
\end{cases}
\label{dase}%
\end{equation}
where $n\geq2$, $\kappa\geq0$ and $\gamma>0$. \ The fractional laplacian
operator $(-\Delta)^{\gamma}$ is defined by
\[
\widehat{\lbrack(-\Delta)^{\gamma}f]}(\xi)=|\xi|^{2\gamma}\widehat{f}(\xi),
\]
where $\widehat{f}=\int_{\mathbb{R}^{n}}e^{-ix\cdot\xi}f(\xi)d\xi$ stands for
the Fourier transform of $f.$ The velocity field $u$ is determined by the
active scalar $\theta$ by means of the multiplier operators%

\begin{equation}
u=P[\theta]=(\widetilde{P}_{1}[\theta],...,\widetilde{P}_{n}[\theta]),
\label{velocity}%
\end{equation}
such that $\nabla\cdot u=0$, and
\begin{equation}
u_{j}=\widetilde{P}_{j}[\theta]=\sum_{i=1}^{n}a_{ij}\mathcal{R}_{i}%
\Lambda^{-1}P_{i}[\theta],\;\;\text{for}\ 1\leq j\leq n, \label{ujotas}%
\end{equation}
\bigskip where $\Lambda=(-\Delta)^{\frac{1}{2}}$, $\mathcal{R}_{i}%
=-\partial_{i}(-\Delta)^{-\frac{1}{2}}$ is the $i$-th Riesz transform,
$a_{ij}$'s are constant and

\begin{equation}
\widehat{P_{i}[\theta]}(\xi)=P_{i}(\xi)\widehat{\theta}(\xi)\text{.}
\label{Pi´s}%
\end{equation}
Denoting $I=\sqrt{-1},$ it follows that
\[
\text{ }\widehat{\widetilde{P}_{j}[\theta]}(\xi)=\widetilde{P}_{j}%
(\xi)\widehat{\theta}(\xi)\text{ with }\widetilde{P}_{j}(\xi)=\sum_{i=1}%
^{n}a_{ij}\frac{\xi_{i}I}{\left\vert \xi\right\vert ^{2}}P_{i}(\xi),
\]
and the vector field $u$ can be expressed in Fourier variables in the shorter
form
\begin{equation}
\widehat{u}=\widehat{P[\theta]}=P(\xi)\widehat{\theta}(\xi)\text{ where }%
P(\xi)=(\widetilde{P}_{1}(\xi),...,\widetilde{P}_{n}(\xi)). \label{p-Fourier}%
\end{equation}
Throughout this manuscript the symbol $P_{i}(\xi)$ in (\ref{Pi´s}) is assumed
to belong to $C^{[\frac{n}{2}]+1}(\mathbb{R}^{n}\backslash\{0\})$ with%

\begin{equation}
\left\vert \frac{\partial^{\alpha}P_{i}}{\partial\xi^{\alpha}}(\xi)\right\vert
\leq C|\xi|^{\beta-|\alpha|}, \label{Pi-cond}%
\end{equation}
for all $\alpha\ \in\ (\mathbb{N}\cup\{0\})^{n}$, $|\alpha|\leq\lbrack\frac
{n}{2}]+1$ and $\xi\neq0$, where $\beta\geq0$. The brackets $[\cdot]$ stands
for the greatest integer function. In particular, for $\alpha=0$ it follows
from (\ref{p-Fourier}) and (\ref{Pi-cond}) that
\begin{equation}
\left\vert \widehat{u}(\xi)\right\vert \leq C\left\vert \xi\right\vert
^{\beta-1}\left\vert \widehat{\theta}(\xi)\right\vert ,\text{ for all }\xi
\neq0. \label{est-field-1}%
\end{equation}
Concerning the criticality of (\ref{dase})-(\ref{ujotas}), there is an
interplay between the field $u$ and fractional viscosity $\left(
-\Delta\right)  ^{\gamma}$ expressed by means of three basic cases:
sub-critical $\beta<2\gamma,$ critical $\beta=2\gamma,$ and super-critical
$\beta>2\gamma.$

We could consider an arbitrary $\kappa>0,$ nevertheless $\kappa=1$ is assumed
for the sake of simplicity. The IVP (\ref{dase})-(\ref{ujotas}) can be
converted into the integral equation
\begin{equation}
\theta(t)=G_{\gamma}(t)\theta_{0}+B(\theta,\theta)(t), \label{mild}%
\end{equation}
where
\begin{equation}
B(\theta,\varphi)(t)=-\int_{0}^{t}G_{\gamma}(t-s)(\nabla_{x}\cdot
(P[\theta]\varphi))(s)ds \label{termo bilinear}%
\end{equation}
and $G_{\gamma}(t)$ is the convolution operator with kernel given in Fourier
variables by $\hat{g}_{\gamma}(\xi,t)=e^{-t|\xi|^{2\gamma}}$. Solutions of
(\ref{mild}) are called mild ones for (\ref{dase})-(\ref{ujotas}).

Assuming that $P_{i}$'s are homogeneous functions of degree $\beta$, we have
formally that
\[
\theta_{\lambda}=\lambda^{2\gamma-\beta}\theta(\lambda x,\lambda^{2\gamma}t)
\]
verifies (\ref{dase})-(\ref{ujotas}), for all$\ \lambda>0,$ provided that
$\theta$ does so. It follows that
\begin{equation}
\theta\rightarrow\theta_{\lambda}=\lambda^{2\gamma-\beta}\theta(\lambda
x,\lambda^{2\gamma}t),\ \text{for}\ \lambda>0, \label{scaling}%
\end{equation}
is the scaling map for (\ref{dase})-(\ref{ujotas}). Also, making
$t\rightarrow0^{+}$ in (\ref{scaling}), one obtains the scaling for the
initial data%

\begin{equation}
\theta_{0}\rightarrow\lambda^{2\gamma-\beta}\theta_{0}(\lambda x).
\label{scaling initial}%
\end{equation}
In view of (\ref{Pi-cond}), even when $P_{i}$ is not homogeneous, we can
consider (\ref{scaling}) as an intrinsic scaling for (\ref{dase}%
)-(\ref{ujotas}) in the sense that it is useful to identify threshold indexes
for functional settings and properties of solutions. One of our aims is to
provide a global well-posedness result for (\ref{dase})-(\ref{ujotas}) in a
scaling invariant framework outside the $L^{2}$-space.

Active scalar equations like (\ref{dase})-(\ref{ujotas}) arise in a large
number of physical models in fluid mechanics and atmospheric science. Examples
of those are 2D surface quasi-geostrophic equation (SQG) $u=\nabla^{\perp
}((-\Delta)^{-1/2}\theta)$ ($\beta=1$), Burgers equation $u=\theta$ ($\beta
=1$)$,$ 2D vorticity equation $u=\nabla^{\perp}(-\Delta)^{-1}\theta$
($\beta=0$). SQG is a famous model with a lot of papers concerning existence,
uniqueness, regularity and asymptotic behavior of solutions in the inviscid
case $\kappa=0$ or in the subcritical ($1/2<\gamma<1$), critical ($\gamma
=1/2$) and supercritical ($\gamma\in(0,1/2)$) ranges. Without making a
complete list, we would like to mention \cite{Caff-Vass}, \cite{Car-Fer1}%
,\cite{Const2},\cite{Const3},\cite{Const4},\cite{Cordoba1},\cite{Dong-1}%
,\cite{Gancedo1},\cite{Ju},\cite{Kiselev1},\cite{Kiselev2},\cite{Kiselev3}%
,\cite{NS},\cite{SS}, and their references. In the case $u=\theta$, see
\cite{Kiselev3} and \cite{Dong-3} for results on blow-up, global existence and
regularity of solutions. One dimensional active scalar models have also
attracted the attention of many authors, see e.g. \cite{Car-Fer3}%
,\cite{Cordoba2},\cite{Dong-2},\cite{Li-Ro} where the reader can find global
existence, finite-time singularity and asymptotic behavior results with
velocity coupled via singular integral operators that are zero-order
multiplier ones.

In the case of SQG, notice that $u$ can be written by using Riesz transform
as
\begin{equation}
u=(-\mathcal{R}_{2}\theta,\mathcal{R}_{1}\theta) \label{Riesz-field}%
\end{equation}
and then the velocity is coupled to the active scalar via zero-order
multiplier operators. The model (\ref{dase})-(\ref{ujotas}) was introduced in
\cite{Chae1},\cite{Chae2} and covers positive-order couplings when $\beta>1$
(see (\ref{est-field-1})). In this range, the operator $P[\cdot]$ behaves
\textquotedblleft morally\textquotedblright\ like a positive derivative of
$(\beta-1)$-order and produces more difficulties in comparison with SQG
($\beta=1$, zero-order) and $\beta<1$ (negative-order).

The paper \cite{Chae1} deals mainly with the inviscid case $\kappa=0$, while
\cite{Chae2} with the dissipative one $\kappa>0$.\ This last work is our main
motivation since we also focus in the dissipative model. The authors of
\cite{Chae2} showed existence of global solutions in $L^{\infty}%
((0,\infty);Y)$ for (\ref{dase})-(\ref{ujotas}) where $Y=L^{1}\cap L^{\infty
}\cap B_{q,\infty}^{s,M}$ with $s>1$ and $2\leq q\leq\infty.$ The index
$M=\{M_{j}\}_{j\geq-1}$ is a sequence and the space $B_{q,\infty}^{s,M}$ is an
extension of the classical Besov space $B_{q,\infty}^{s}$ where the
$B_{q,\infty}^{s,M}$-norm increases according to the growth of $M$. The
results of \cite{Chae2} consider couplings $P[\cdot]$ in (\ref{velocity}) such
that $P_{i}\in C^{\infty}(\mathbb{R}^{n}\backslash\{0\}),$ $P_{i}$ is radially
symmetric, $P_{i}=P_{i}(\left\vert \xi\right\vert )$ is nondecreasing with
$\left\vert \xi\right\vert $, and a technical growth hypothesis involving
$P_{i}(\xi)$ and the sequence $M$. Applying their results to the special case
\begin{equation}
u=\nabla^{\perp}(\Lambda^{\beta-2}\theta)=\Lambda^{\beta-1}(-\mathcal{R}%
_{2}\theta,\mathcal{R}_{1}\theta) \label{field1}%
\end{equation}
with $n=2$, $0\leq\beta<2\gamma<1$ (within the sub-critical range), $\kappa>0$
and $M_{j}=j+1$, they obtained well-posedness of solutions with initial data
in $L^{1}\cap L^{\infty}\cap B_{q,\infty}^{s,M}.$ Roughly speaking, the
technique employed in \cite{Chae2} for constructing solutions relies
on\textit{ }a successive approximation scheme together \textit{a priori}
estimates involving Besov norms. The field (\ref{field1}) corresponds to the
modified SQG that interpolates 2D vorticity equation and SQG by varying the
parameter $\beta$ from $0$ to $1.$ This model has been studied for instance in
\cite{Chae1},\cite{Const1},\cite{Kiselev3},\cite{May1},\cite{Miao1}%
,\cite{Miao-2} where one can find existence and regularity results with data
in Sobolev spaces $H^{m}$ with $m\geq0$. The conditions $\kappa>0$, $\beta
\in\lbrack0,1]$ and $\beta=2\gamma$ were assumed in \cite{Const1}%
,\cite{Kiselev3},\cite{May1},\cite{Miao1}; $\kappa>0$ and $1\leq\beta
<2\gamma<2$ in \cite{Miao-2}; and $\kappa=0$ and $\beta\in\lbrack1,2]$ in
\cite{Chae1}. In this last work, local well-posedness of $H^{m}(\mathbb{R}%
^{2})$-solutions was proved for (\ref{dase})-(\ref{field1}) with $m\geq4$.

In this paper we prove the global-in-time well-posedness of (\ref{dase}%
)-(\ref{ujotas}) in the Lebesgue space $L^{\frac{n}{2\gamma-\beta}}%
(\mathbb{R}^{n})$ without smallness conditions (see Theorem \ref{teoglobal}).
This is the unique $L^{r}$-space whose norm is invariant by the scaling
(\ref{scaling initial}), that is, $L^{\frac{n}{2\gamma-\beta}}$ is the
critical one in the scale of Lebesgue spaces. We can consider initial data
outside the $L^{2}$-framework and, due to the inclusion $L^{1}\cap L^{\infty
}\subset L^{\frac{n}{2\gamma-\beta}}$, our initial data class is larger than
that of \cite{Chae2}. In comparison with \cite{Chae2}, some new symbols
$P_{i}(\xi)$ are considered here (e.g. non-radially symmetric ones). Even for
a singular initial data $\theta_{0}\in L^{\frac{n}{2\gamma-\beta}}%
(\mathbb{R}^{n})$, the global solution $\theta\in BC([0,\infty);L^{\frac
{n}{2\gamma-\beta}}(\mathbb{R}^{n}))$ is instantaneously $C^{\infty}$-smoothed
out and verifies (\ref{dase})-(\ref{ujotas}) classically, for all $t>0$. Here
we focus in the range $\beta\geq1$ and consider the sub-critical case
$\beta<2\gamma$. More precisely, we assume
\begin{equation}
1\leq2\beta-1<2\gamma<\min\{\frac{2}{3}(n+\beta+1),(n+1)\}. \label{cond-1}%
\end{equation}
The range $0\leq\beta<2\gamma$ with $\beta<1$ also can be treated with an
adaptation on the proofs (see Remark \ref{rem1}).

Also, we show some decay properties in $L^{q}$-norms (see Theorem
\ref{teoglobal}). Precisely, for $\frac{n}{2\gamma-\beta}\leq q\leq\infty$ and
$\theta_{0}\in L^{\frac{n}{2\gamma-\beta}}(\mathbb{R}^{n}),$ we obtain the
time-polynomial decay
\begin{equation}
\left\Vert \theta(\cdot,t)\right\Vert _{L^{q}}\leq Ct^{-(\frac{2\gamma-\beta
}{2\gamma}-\frac{n}{2\gamma q})}\,,\ \text{for all }t>0\text{.} \label{decay1}%
\end{equation}
Assuming further $\theta_{0}\in L^{\frac{n}{2\gamma-\beta}}(\mathbb{R}%
^{n})\cap L^{1}(\mathbb{R}^{n}),$ the solution $\theta$ belongs to
$BC([0,\infty);L^{1}(\mathbb{R}^{n}))$ and the estimate (\ref{decay1}) is
improved to%
\begin{equation}
\left\Vert \theta(\cdot,t)\right\Vert _{L^{q}}\leq Ct^{-(\frac{2\gamma-\beta
}{2\gamma}-\frac{n}{2\gamma q})-(\frac{n+\beta}{2\gamma}-1)}\,,\text{ for all
}t>0, \label{decay2}%
\end{equation}
where $1\leq q\leq\infty.$ Notice that the decay in (\ref{decay2}) is faster
than those of (\ref{decay1}) due to the condition $2\gamma<\frac{2}{3}%
(n+\beta+1)<n+\beta$.

In view of the $L^{p}$-$L^{q}$ estimate (\ref{est linear2}) for the semigroup
$G_{\gamma}(t)$, it is not expected that (\ref{decay1}) holds true for
$q<\frac{n}{2\gamma-\beta}$ and an arbitrary $\theta_{0}\in L^{\frac
{n}{2\gamma-\beta}}(\mathbb{R}^{n}).$ Thus $\theta(\cdot,t)$ may not be a
$L^{2}$-solution when $2<\frac{n}{2\gamma-\beta}$ although $\theta(\cdot,t)\in
C^{\infty}(\mathbb{R}^{n})$, for all $t>0.$ Even in the subcritical case, this
fact seems to prevent an adaptation from previous techniques based on $L^{2}%
$-frameworks  (see e.g. the famous papers \cite{Caff-Vass, Kiselev1}) in order
to obtain global well-posedness of $L^{\frac{n}{2\gamma-\beta}}(\mathbb{R}%
^{n})$-solutions. Roughly speaking, the approach employed here relies on
time-weighted Kato type norms, scaling arguments, and arguments of the type
parabolic De Giorgi-Nash-Moser. These ingredients also were used in
\cite{Car-Fer1} in order to analyze SQG ($\beta=1$). However, due to the
coupling between $\theta$ and $u$ being via a positive-order operator, the
model (\ref{dase})-(\ref{ujotas}) requires more involved arguments and further
care in comparison with SQG. For instance, since $P[\cdot]$ is not continuous
from $L^{p_{1}}$ to $L^{p_{2}}$ when $\beta>1$, we need to employ an auxiliary
time-weighted Kato-type norm based on homogeneous Sobolev spaces $\dot{H}%
_{q}^{s}$ with $q>\frac{n}{2\gamma-\beta}$ in order to control the nonlinear
term in (\ref{dase})-(\ref{ujotas}). So, different from SQG, Sobolev norms
play here a crucial role for the local existence and extension of solutions
with data in Lebesgue spaces (see e.g. (\ref{aux-prop1}) and (\ref{aux-def1}%
)-(\ref{aux-ext-10}), respectively). Let us also mention that there is no
maximum principle for $\dot{H}_{q}^{s}$-norms when $s>0$; and consequently
there is a lack of global-in-time control on these norms (see
(\ref{Contituity-Hq})).

In view of the structure of (\ref{dase}), it is natural to wonder about
symmetry properties of solutions under symmetry conditions for the symbols
$P_{i}(\xi)$ and initial data $\theta_{0}.$ In Theorem \ref{Teo-sym}, we show
that the global solution given in Theorem \ref{teoglobal} is radially
symmetric, for all $t>0$, provided that $\theta_{0}$ and $div_{\xi}(P(\xi))$
present this same property. Moreover, results on odd and even symmetry of
solutions are obtained under parity conditions for $\theta_{0}$ and $P_{i}$'s.
In Remark \ref{rem2}, we also comment about conditions for solutions to be non-symmetric.

Let us also comment on \textit{log}-type couplings which are interesting ones
covered by (\ref{dase})-(\ref{ujotas}). Ohkitani \cite{Okitani} has presented
numerical evidences that, even with $\kappa=0,$ (\ref{dase}) with $n=2$ and
\begin{equation}
u=\nabla^{\perp}(\log(I-\Delta))^{\chi}\theta,\text{ }\chi>0,
\label{log-field-1}%
\end{equation}
may be globally well-posed. The authors of \cite{Chae1} have proved local
well-posedness of $H^{4}(\mathbb{R}^{2})$-solutions for (\ref{dase}%
)-(\ref{log-field-1}) with $\kappa>0$. As pointed in \cite{Chae1}, the field
(\ref{log-field-1}) is of order higher (at least logarithmically) than
derivatives of order $1$ and in particular than (\ref{Riesz-field}). Another
examples are
\begin{align}
P_{i}(\xi)  &  =\left\vert \xi\right\vert ^{\sigma}(\log(1+\left\vert
\xi\right\vert ^{2}))^{\chi},\text{ }\chi\geq0,\text{ }\label{log-field-2}\\
P_{i}(\xi)  &  =\left\vert \xi\right\vert ^{\sigma}(\log(1+\log(1+\left\vert
\xi\right\vert ^{2})))^{\chi},\chi\geq0, \label{log-field-3}%
\end{align}
which are indeed of order higher than (\ref{log-field-1}) when $\sigma>1$.
These couplings are also treated in \cite{Chae2} with $\sigma=\beta$ and
$\chi\geq0.$ When $\sigma=0$ and $n=2,$ (\ref{log-field-2}) and
(\ref{log-field-3}) correspond to \textit{log} and \textit{log-log}
Navier-Stokes which are intermediate models between 2D vorticity equation and
SQG. See \cite{Chae3} for further details and global existence results in the
case $\kappa=0,$ $0\leq\chi\leq1$ and data $\theta_{0}\in L^{1}\cap L^{\infty
}\cap B_{q,\infty}^{s}$, where $B_{q,\infty}^{s}$ stands for an inhomogeneous
Besov space with $s>1$ and $q>2$. An interest in\ \textit{log}-type couplings
has also arisen in connection with other fluid mechanics models (see
\cite{Chae5}).

Finally, we remark that our results cover the couplings (\ref{field1}),
(\ref{log-field-1}), (\ref{log-field-2}) and (\ref{log-field-3}). The
condition (\ref{Pi-cond}) is clearly satisfied by (\ref{field1}), and if
$\beta\in\lbrack1,2]$ and $2\beta-1<2\gamma<\min\{2+\frac{2\beta}{3},3\}$ then
(\ref{cond-1}) holds true. Also, (\ref{log-field-1}) verifies (\ref{Pi-cond})
with $\beta=1+\varepsilon,$ for any $\varepsilon>0,$ and we have
(\ref{cond-1}) when $\frac{1}{2}<\gamma<\frac{4}{3}$ and $0<\varepsilon
<\gamma-\frac{1}{2}$. By considering $\beta=\sigma+\varepsilon,$ conditions
analogous to the ones for (\ref{log-field-1}) can be obtained for
(\ref{log-field-2}) and (\ref{log-field-3}) with $\chi>0.$ The cases
(\ref{log-field-2}) and (\ref{log-field-3}) with $\chi=0$ are similar to
(\ref{field1}).

This manuscript is organized as follows. In the next section we recall some
estimates in $L^{q}(\mathbb{R}^{n})$ and Sobolev homogeneous spaces for
Fourier multiplier operators and the semigroup $\{G_{\gamma}(t)\}_{t\geq0}$.
Our results are stated in section 3 in two theorems, namely Theorems
\ref{teoglobal} and \ref{Teo-sym}. Estimates for the bilinear operator
(\ref{termo bilinear}) are obtained in section 4. Local well-posedness and
some properties of solutions are proved in subsection 5.1. The proofs of
Theorems \ref{teoglobal} and \ref{Teo-sym} are performed in subsections 5.2
and 5.3, respectively.

\section{Preliminaries}

\hspace{0.55cm} In this section we recall some estimates for the fundamental
solution of the linear part of (\ref{dase}) in $L^{p}(\mathbb{R}^{n})$ and
$\dot{H}_{p}^{s}(\mathbb{R}^{n})$, whose norms will be denoted by $\Vert
\cdot\Vert_{p}$ and $\Vert\cdot\Vert_{\dot{H}_{p}^{s}}$, respectively.

We remember that given $s\in\mathbb{R}$ and $1<p<\infty$, the homogeneous
Sobolev space $\dot{H}_{p}^{s}(\mathbb{R}^{n})$ is the space of all
$u\in\mathcal{S}^{\prime}/\mathcal{P}$ such that $(-\Delta)^{\frac{s}{2}}u\in$
$L^{p}(\mathbb{R}^{n})$. In other words, $\dot{H}_{p}^{s}=(-\Delta)^{-\frac
{s}{2}}L^{p}$ and it is a Banach space with norm
\[
\Vert u\Vert_{\dot{H}_{p}^{s}}=\Vert(-\Delta)^{\frac{s}{2}}u\Vert_{p}.
\]
The following Sobolev type embedding holds true%

\begin{equation}
\dot{H}_{p_{2}}^{s_{2}}(\mathbb{R}^{n})\subset\dot{H}_{p_{1}}^{s_{1}%
}(\mathbb{R}^{n}), \label{Sobolev}%
\end{equation}
for $1<p_{2}\leq p_{1}<\infty$ and $s_{1}-\frac{n}{p_{1}}=s_{2}-\frac{n}%
{p_{2}}$. The reader is refereed to \cite[chapter 6]{Grafakos} for further
details on these spaces.

The next lemma gives estimates for certain multiplier operators acting in
$\dot{H}_{p}^{s}(\mathbb{R}^{n})$ (see e.g. \cite{Kozo1})

\begin{lemma}
\label{multiplier} Let $m$, $s\in\mathbb{R}$, $1<p<\infty$ and $F(\xi)\in
C^{[\frac{n}{2}]+1}(\mathbb{R}^{n}\backslash\{0\})$, where $[\cdot]$ stands
for the greatest integer function. Assume that there is $L>0$ such that
\begin{equation}
\left\vert \frac{\partial^{\alpha}F}{\partial\xi^{\alpha}}(\xi)\right\vert
\leq L|\xi|^{m-|\alpha|}, \label{est fjotas}%
\end{equation}
for all $\alpha\in(\mathbb{N}\cup\{0\})^{n}$, $|\alpha|\leq\lbrack\frac{n}%
{2}]+1$, and $\xi\neq0$. Then the multiplier operator $F(D)$ on $\mathcal{S}%
^{\prime}/\mathcal{P}$ is bounded from $\dot{H}_{p}^{s}$ to $\dot{H}_{p}%
^{s-m}$. Moreover, the following estimate holds true
\begin{equation}
\Vert F(D)f\Vert_{\dot{H}_{p}^{s-m}}\leq C\Vert f\Vert_{\dot{H}_{p}^{s}},
\label{fjotas}%
\end{equation}
where $C>0$ is independent of $f$.
\end{lemma}

The next lemma gives estimates for $\{G_{\gamma}(t)\}_{t\geq0}$ on spaces
$L^{p}(\mathbb{R}^{n})$ and $\dot{H}_{p}^{s}(\mathbb{R}^{n})$.

\begin{lemma}
\label{sol fund} Let $n\geq2,$ $0<\gamma<\infty$, $1\leq p\leq q\leq\infty$
and $k\in(\mathbb{N}\cup\{0\})^{n}.$ Then
\begin{equation}
\Vert\nabla_{x}^{k}G_{\gamma}(t)f\Vert_{q}\leq C\ t^{-\frac{\left\vert
k\right\vert }{2\gamma}-\frac{n}{2\gamma}(\frac{1}{p}-\frac{1}{q})}\Vert
f\Vert_{p}\text{,} \label{est linear2}%
\end{equation}
for all $f\in L^{p}(\mathbb{R}^{n})$. Now, let $s_{1}\leq s_{2}$, $s_{i}%
\in\mathbb{R}$ and $1<p_{1}\leq p_{2}<\infty$. There is a constant $C>0$ such
that
\begin{equation}
\Vert G_{\gamma}(t)f\Vert_{\dot{H}_{p_{2}}^{s_{2}}}\leq Ct^{-\frac
{(s_{2}-s_{1})}{2\gamma}-\frac{n}{2\gamma}(\frac{1}{p_{1}}-\frac{1}{p_{2}}%
)}\Vert f\Vert_{\dot{H}_{p_{1}}^{s_{1}}}. \label{est linear}%
\end{equation}
for all $f\in\dot{H}_{p_{1}}^{s_{1}}$. Moreover, given $f\in L^{\frac
{n}{2\gamma-\beta}}(\mathbb{R}^{n})$ with $1\leq\beta<2\gamma\leq n+\beta$ and
$\frac{n}{2\gamma-\beta}<q<\infty$, then
\begin{equation}
\sup_{0<t<T}t^{\eta_{q}}\Vert G_{\gamma}(t)f\Vert_{\dot{H}_{q}^{\beta-1}}\leq
C\Vert f\Vert_{\frac{n}{2\gamma-\beta}}\qquad\mbox{and}\qquad\lim
_{t\rightarrow0^{+}}t^{\eta_{q}}\Vert G_{\gamma}(t)\theta_{0}\Vert_{\dot
{H}_{q}^{\beta-1}}=0 \label{zeroPrincipal}%
\end{equation}
where $\eta_{q}=\frac{2\gamma-1}{2\gamma}-\frac{n}{2\gamma q}$ and $C$ is
independent of $f$ and $0<T\leq\infty$.
\end{lemma}

\textbf{Proof. }The estimate (\ref{est linear2}) is well-known (see e.g.
\cite{Car-Fer1} for a proof). Also, (\ref{est linear2}) still holds true by
replacing $\nabla_{x}^{k}$ by $(-\Delta)^{\frac{\left\vert k\right\vert }{2}%
}.$ In view of the latter comment and $(-\Delta)^{\frac{s_{2}}{2}}%
=(-\Delta)^{\frac{s_{2}-s_{1}}{2}}(-\Delta)^{\frac{s_{1}}{2}}$, we obtain
(\ref{est linear}) from (\ref{est linear2}) because $G_{\gamma}(t)$ commutates
with $(-\Delta)^{\frac{s_{1}}{2}}.$ The estimate in (\ref{zeroPrincipal})
comes from (\ref{est linear}) with $p_{2}=q,$ $s_{2}=\beta-1,$ $p_{1}=\frac
{n}{2\gamma-\beta}$ and $s_{1}=0.$ Due to (\ref{est linear}), it is easy to
see that the limit in (\ref{zeroPrincipal}) holds true for every $\theta
_{0}\in L^{\frac{n}{2\gamma-\beta}}\cap\dot{H}_{q}^{\beta-1}.$ This fact
together with $\overline{L^{\frac{n}{2\gamma-\beta}}\cap\dot{H}_{q}^{\beta-1}%
}^{\left\Vert \cdot\right\Vert _{\frac{n}{2\gamma-\beta}}}=L^{\frac{n}%
{2\gamma-\beta}}$ and the estimate in (\ref{zeroPrincipal}) yield the limit in
(\ref{zeroPrincipal}), for every $\theta_{0}\in L^{\frac{n}{2\gamma-\beta}%
}(\mathbb{R}^{n}).$

{\vskip 0pt \hfill\hbox{\vrule height 5pt width 5pt depth 0pt} \vskip 12pt}

\section{Results}

\hspace{0.55cm} This section is devoted to the statements of the results whose
proofs will be performed in section 5.

\begin{theorem}
\label{teoglobal} (Global solutions) Assume the condition (\ref{cond-1}) and
let $\eta_{q}=\frac{2\gamma-1}{2\gamma}-\frac{n}{2\gamma q}$ and $\tilde{\eta
}_{q}=\frac{2\gamma-\beta}{2\gamma}-\frac{n}{2\gamma q}.$ If $\theta_{0}\in
L^{\frac{n}{2\gamma-\beta}}(\mathbb{R}^{n})$ then there is a unique global
solution $\theta\in BC([0,\infty);L^{\frac{n}{2\gamma-\beta}}(\mathbb{R}%
^{n}))$ for (\ref{dase})-(\ref{ujotas}) such that
\begin{align}
t^{\tilde{\eta}_{q}}\theta &  \in BC\left(  (0,\infty),L^{q}(\mathbb{R}%
^{n})\right)  ,\text{ for all }\frac{n}{2\gamma-\beta}<q\leq\infty
,\label{Decay-Lq-1}\\
t^{\eta_{q}}\theta &  \in C((0,\infty);\dot{H}_{q}^{\beta-1}(\mathbb{R}%
^{n})),\text{ \ for all }\frac{n}{2\gamma-\beta}<q<\infty,
\label{Contituity-Hq}%
\end{align}
where the limits of $t^{\tilde{\eta}_{q}}\theta$ in (\ref{Decay-Lq-1}) and
$t^{\eta_{q}}\theta$ in (\ref{Contituity-Hq}) are zero as $t\rightarrow0^{+}.$
Moreover, if $\theta_{0}\in L^{1}(\mathbb{R}^{n})\cap L^{\frac{n}%
{2\gamma-\beta}}(\mathbb{R}^{n})$ and $1<q\leq\infty$, then $\theta
\in\ BC([0,\infty);L^{1}(\mathbb{R}^{n})\cap L^{\frac{n}{2\gamma-\beta}%
}(\mathbb{R}^{n}))$ and
\begin{equation}
t^{\tilde{\eta}_{q}+\frac{n+\beta}{2\gamma}-1}\theta\in BC((0,\infty
);L^{q}(\mathbb{R}^{n})). \label{Decay-Lq-2}%
\end{equation}

\end{theorem}

\begin{remark}
\label{ContDepend}(Continuous dependence on initial data) The proof of Theorem
\ref{teoglobal} also gives that the solution $\theta$ depends continuously on
the data $\theta_{0}$ in finite time intervals $[0,T]$. More precisely, if
$\theta_{k,0}\rightarrow\theta_{0}$ in $L^{\frac{n}{2\gamma-\beta}}%
(\mathbb{R}^{n})$ then $\theta_{k}\rightarrow\theta$ in $C([0,T];L^{\frac
{n}{2\gamma-\beta}}(\mathbb{R}^{n}))$, for all $0<T<\infty$, where $\theta
_{k}$ is the solution with initial data $\theta_{k,0}.$
\end{remark}

\begin{remark}
\label{rem1}One can treat the range $0\leq\beta<1$ by modifying the estimates
of section 4 (particularly (\ref{bili-1}) and (\ref{bili-4})). For that
matter, one should replace $\sup_{0<t<T}t^{\eta_{l}}\left\Vert \theta
(\cdot,t)\right\Vert _{\dot{H}_{l}^{\beta-1}}$ by $\ \sup_{0<t<T}%
t^{\tilde{\eta}_{l}}\left\Vert \theta(\cdot,t)\right\Vert _{l}$ into those
estimates $(l=q,r)$. In fact, due to Hardy-Littlewood-Sobolev inequality, this
case is easier to handling than $\beta>1$ and it is not necessary to use norms
based on homogeneous Sobolev spaces in order to prove existence of solutions.
\end{remark}

Before proceeding, we recall the concept of even and odd symmetry. A function
$h$ is said to be even (resp. odd) when $h(x)=h(-x)$ (resp. $h(x)=-h(-x)$).

\begin{theorem}
\label{Teo-sym}(Symmetry) Under the hypotheses of Theorem \ref{teoglobal}.

\begin{description}
\item[(i)] \textit{The solution }$\theta(x,t)$\textit{\ is odd (resp. even)
for all }$t>0,$\textit{ }provided that $\theta_{0}$ and $P_{i}$'s are odd
(resp. even).

\item[(ii)] Let $P(\xi)$ be as in (\ref{p-Fourier}). If $\theta_{0}$ and
$div_{\xi}(P(\xi))$ \textit{are radially symmetric then }$\theta
(x,t)$\textit{\ is radially symmetric for all }$t>0$.
\end{description}
\end{theorem}

\bigskip

\begin{remark}
\label{rem2}\textit{(non-symmetry)} Adapting the arguments in the proof of
Theorem \ref{Teo-sym}, one also can prove the following non-symmetry results:
if $\theta_{0}$ is odd (resp. even) and $P_{i}$'s are even (resp. odd) then
$\theta(x,t)$ is not odd (resp. not even). Also, if $\theta_{0}$ is nonradial
and $div_{\xi}(P(\xi))$ is radial, then $\theta(x,t)$ is not radially
symmetric. The detailed verification is left to the reader.
\end{remark}

\section{\bigskip Bilinear Estimates}

\hspace{0.55cm} This part of the article is devoted to estimates for the
bilinear term (\ref{termo bilinear}).

\begin{lemma}
\label{LemaBiliq} Let $0<T\leq\infty,$ $n\geq2,$ $1\leq\beta<2\gamma<\infty,$
and let $1<q<\infty$ be such that $\frac{\beta-1}{n}<\frac{1}{q}<\frac
{2\gamma-1}{n}.$ Denote $\eta_{q}=\frac{2\gamma-1}{2\gamma}-\frac{n}{2\gamma
q}$ and $\tilde{\eta}_{q}=\frac{2\gamma-\beta}{2\gamma}-\frac{n}{2\gamma q}$.

\begin{itemize}
\item[(i)] If $\frac{2\gamma-(\beta+1)}{n}-\frac{1}{q}<\frac{1}{r}\leq\frac
{1}{q^{\prime}}$ and $q^{\prime}\leq p\leq\infty$ then there are positive
constants $K_{1},K_{2},K_{3},$ independent of $\theta,\phi$ and $T,$ such
that
\begin{align}
\sup_{0<t<T}t^{\tilde{\eta}_{r}}\Vert B(\theta,\phi)\Vert_{r}\text{ }  &  \leq
K_{1}\sup_{0<t<T}t^{\eta_{q}}\Vert\theta\Vert_{\dot{H}_{q}^{\beta-1}}%
\,\sup_{0<t<T}t^{\tilde{\eta}_{r}}\Vert\phi\Vert_{r},\label{bili-1}\\
\sup_{0<t<T}\Vert B(\theta,\phi)\Vert_{p}\text{ }  &  \leq K_{2}\sup
_{0<t<T}t^{\eta_{q}}\Vert\theta\Vert_{\dot{H}_{q}^{\beta-1}}\,\sup
_{0<t<T}\Vert\phi\Vert_{p},\label{bili-2}\\
\sup_{0<t<T}\Vert B(\theta,\phi)\Vert_{1}  &  \leq K_{3}T^{\frac{2\gamma
-1}{2\gamma}-\eta_{q}}\sup_{0<t<T}t^{\eta_{q}}\Vert\theta\Vert_{\dot{H}%
_{q}^{\beta-1}}\,\sup_{0<t<T}\Vert\phi\Vert_{q^{\prime}}. \label{bili-3}%
\end{align}

\item[(ii)] If $\frac{2\gamma-2}{n}-\frac{1}{q}<\frac{1}{r}<\frac{n+\beta
-1}{n}-\frac{1}{q}$ then
\begin{equation}
\sup_{0<t<T}t^{\eta_{r}}\Vert B(\theta,\phi)\Vert_{\dot{H}_{r}^{\beta-1}%
}\text{ }\leq K_{4}\sup_{0<t<T}t^{\eta_{r}}\Vert\theta\Vert_{\dot{H}%
_{r}^{\beta-1}}\,\sup_{0<t<T}t^{\eta_{q}}\Vert\phi\Vert_{\dot{H}_{q}^{\beta
-1}}, \label{bili-4}%
\end{equation}
where $K_{4}>0$ is a constant independent of $\theta,\phi$ and $T.$
\end{itemize}
\end{lemma}

\bigskip

\textbf{Proof. }

\bigskip

\textit{Proof of part (i):} Let $p_{1}=p$ and $p_{2}=r$. Using Lemma
\ref{sol fund} and H\"{o}lder inequality, we estimate
\begin{align}
\Vert B(\theta,\phi)\Vert_{p_{i}}  &  \leq\int_{0}^{t}\Vert\nabla_{x}%
G_{\gamma}(t-s)(P[\theta]\phi)(s)\Vert_{p_{i}}\text{ }ds\nonumber\\
&  \leq C\int_{0}^{t}(t-s)^{-\frac{1}{2\gamma}-\frac{n}{2\gamma q}}%
\,\Vert(P[\theta]\phi)(s)\Vert_{\frac{p_{i}q}{p_{i}+q}}\text{ }ds\nonumber\\
&  \leq C\int_{0}^{t}(t-s)^{-\frac{1}{2\gamma}-\frac{n}{2\gamma q}}\,\Vert
P[\theta(s)]\Vert_{q}\,\Vert\phi(s)\Vert_{p_{i}}\text{ }ds\nonumber\\
&  \leq C\int_{0}^{t}(t-s)^{-\frac{1}{2\gamma}-\frac{n}{2\gamma q}}%
\,\Vert\theta(s)\Vert_{\dot{H}_{q}^{\beta-1}}\,\Vert\phi(s)\Vert_{p_{i}}\text{
}ds \label{Bili}%
\end{align}
where $i=1,2$ and in the third line we have used Lemma \ref{multiplier} in
order to infer
\[
\Vert P[\theta]\Vert_{q}\leq C\Vert\theta\Vert_{\dot{H}_{q}^{\beta-1}}.
\]
Therefore
\begin{align}
\Vert B(\theta,\phi)\Vert_{p}  &  \leq C\,I_{1}(t)\,\sup_{0<t<T}\Vert
\phi(t)\Vert_{p}\,\sup_{0<t<T}t^{\eta_{q}}\Vert\theta(t)\Vert_{\dot{H}%
_{q}^{\beta-1}},\label{aux-bili-1}\\
\Vert B(\theta,\phi)\Vert_{r}  &  \leq C\,I_{2}(t)\,\sup_{0<t<T}t^{\tilde
{\eta}_{r}}\Vert\phi(t)\Vert_{r}\sup_{0<t<T}t^{\eta_{q}}\Vert\theta
(t)\Vert_{\dot{H}_{q}^{\beta-1}}, \label{aux-bili-2}%
\end{align}
where the integrals $I_{1}(t)$ and $I_{2}(t)$ can be computed as
\begin{align}
I_{1}(t)  &  =\int_{0}^{t}(t-s)^{-\frac{1}{2\gamma}-\frac{n}{2\gamma q}%
}s^{-\eta_{q}}\text{ }ds=\int_{0}^{1}(1-s)^{\eta_{q}-1}s^{-\eta_{q}}\text{
}ds=C<\infty,\label{aux-int-1}\\
I_{2}(t)  &  =\int_{0}^{t}(t-s)^{-\frac{1}{2\gamma}-\frac{n}{2\gamma q}%
}s^{-\eta_{q}-\tilde{\eta}_{r}}\,ds=t^{\eta_{q}-1-\eta_{q}-\tilde{\eta}_{r}%
+1}\int_{0}^{1}(1-s)^{\eta_{q}-1}s^{-\eta_{q}-\tilde{\eta}_{r}}%
\,ds\,=\,C\,t^{-\tilde{\eta}_{r}}. \label{aux-int-2}%
\end{align}
The estimates (\ref{bili-1}) and (\ref{bili-2}) follows from (\ref{aux-bili-2}%
) with (\ref{aux-int-2}), and (\ref{aux-bili-1}) with (\ref{aux-int-1}),
respectively.\bigskip\ 

Moreover, we have that
\begin{align*}
\Vert B(\theta,\phi)(t)\Vert_{1}  &  \leq\int_{0}^{t}\Vert\nabla_{x}G_{\gamma
}(t-s)(P[\theta]\phi)(s)\Vert_{1}\text{ }ds\\
&  \leq C\int_{0}^{t}(t-s)^{-\frac{1}{2\gamma}}\,\Vert P[\theta]\Vert
_{q}\,\Vert\phi\Vert_{q^{\prime}}ds\\
&  \leq C\int_{0}^{t}(t-s)^{-\frac{1}{2\gamma}}\,\Vert\theta\Vert_{\dot{H}%
_{q}^{\beta-1}}\,\Vert\phi\Vert_{q^{\prime}}ds\\
&  \leq C\int_{0}^{t}(t-s)^{-\frac{1}{2\gamma}}\,s^{-\eta_{q}}ds\sup
_{0<t<T}t^{\eta_{q}}\Vert\theta\Vert_{\dot{H}_{q}^{\beta-1}}\,\sup
_{0<t<T}\Vert\phi\Vert_{q^{\prime}}\\
&  \leq CT^{1-\frac{1}{2\gamma}-\eta_{q}}\sup_{0<t<T}t^{\eta_{q}}\Vert
\theta\Vert_{\dot{H}_{q}^{\beta-1}}\,\sup_{0<t<T}\Vert\phi\Vert_{q^{\prime}},
\end{align*}
which gives (\ref{bili-3}).

\textit{Proof of part (ii):} Consider
\begin{equation}
\frac{1}{h}=\frac{1}{r}+\frac{1}{q}-\frac{\beta-1}{n}\text{ and }\delta
=\frac{n}{h}-\frac{n}{r}. \label{aux-param1}%
\end{equation}
Note that $\frac{\beta+\delta}{2\gamma}<1$ because $\frac{1}{q}<\frac
{2\gamma-1}{n}$. We employ the continuous inclusion $\dot{H}_{h}%
^{\beta-1+\delta}\subset\dot{H}_{r}^{\beta-1}$, (\ref{est linear}) and
afterwards (\ref{fjotas}) to obtain
\begin{align}
\Vert B(\theta,\phi)\Vert_{\dot{H}_{r}^{\beta-1}}  &  \leq\int_{0}^{t}\Vert
G_{\gamma}(t-s)[\nabla\cdot(P[\theta]\phi)(s)]\Vert_{\dot{H}_{r}^{\beta-1}%
}\text{ }ds\nonumber\\
&  \leq\int_{0}^{t}\Vert G_{\gamma}(t-s)[\nabla\cdot(P[\theta]\phi
)(s)]\Vert_{\dot{H}_{h}^{\beta-1+\delta}}\text{ }ds\nonumber\\
&  \leq\int_{0}^{t}(t-s)^{-\frac{\beta+\delta}{2\gamma}}\Vert\nabla
\cdot(P[\theta]\phi)(s)]\Vert_{\dot{H}_{h}^{-1}}\text{ }ds\nonumber\\
&  \leq\int_{0}^{t}(t-s)^{-\frac{\beta+\delta}{2\gamma}}\Vert P[\theta
]\phi(s)\Vert_{h}\text{ }ds. \label{bilinormh3}%
\end{align}
In view of (\ref{aux-param1}), we can choose $1<l<\infty$ in such a way that
$l>q$, $\frac{1}{h}=\frac{1}{r}+\frac{1}{l}$ and $\frac{1}{l}=\frac{1}%
{q}-\frac{\beta-1}{n}$. Then, H\"{o}lder inequality, (\ref{fjotas}), and
Sobolev embedding (\ref{Sobolev}) imply that
\begin{align}
\Vert P[\theta]\phi\Vert_{h}  &  \leq\Vert P[\theta]\Vert_{r}\Vert\phi
\Vert_{l}\nonumber\\
&  \leq\Vert\theta\Vert_{\dot{H}_{r}^{\beta-1}}\Vert\phi\Vert_{\dot{H}%
_{q}^{\beta-1}}. \label{holdersobolev}%
\end{align}
Inserting (\ref{holdersobolev}) into (\ref{bilinormh3}), we get
\begin{align*}
\Vert B(\theta,\phi)\Vert_{\dot{H}_{r}^{\beta-1}}  &  \leq C\int_{0}%
^{t}(t-s)^{-\frac{\beta+\delta}{2\gamma}}\Vert\theta\Vert_{\dot{H}_{r}%
^{\beta-1}}\Vert\phi\Vert_{\dot{H}_{q}^{\beta-1}}\text{ }ds\\
&  \leq C\int_{0}^{t}(t-s)^{-\frac{\beta+\delta}{2\gamma}}\text{ }s^{-\eta
_{r}-\eta_{q}}ds\sup_{0<t<T}t^{\eta_{r}}\Vert\theta(t)\Vert_{\dot{H}%
_{r}^{\beta-1}}\sup_{0<t<T}t^{\eta_{q}}\Vert\phi(t)\Vert_{\dot{H}_{q}%
^{\beta-1}}\\
&  \leq t^{-\frac{\beta+\delta}{2\gamma}-\eta_{r}-\eta_{q}+1}\int_{0}%
^{1}(1-s)^{-\frac{\beta+\delta}{2\gamma}}s^{-\eta_{r}-\eta_{q}}\text{ }%
ds\sup_{0<t<T}t^{\eta_{r}}\Vert\theta(t)\Vert_{\dot{H}_{r}^{\beta-1}}%
\sup_{0<t<T}t^{\eta_{q}}\Vert\phi(t)\Vert_{\dot{H}_{q}^{\beta-1}}\\
&  \leq Ct^{-\eta_{r}}\sup_{0<t<T}t^{\eta_{r}}\Vert\theta(t)\Vert_{\dot{H}%
_{r}^{\beta-1}}\sup_{0<t<T}t^{\eta_{q}}\Vert\phi(t)\Vert_{\dot{H}_{q}%
^{\beta-1}},
\end{align*}
which is equivalent to (\ref{bili-4}). {\vskip0pt \hfill
\hbox{\vrule height 5pt width 5pt depth 0pt} \vskip12pt}

\section{Proofs}

\subsection{Local in Time Solutions}

\hspace{0.55cm}We start by recalling an abstract lemma in Banach spaces which
is useful in order to avoid extensive fixed point computations (see e.g.
\cite[Theorem 9]{Lewis}).

\begin{lemma}
\label{genlem} Let $X$ be a Banach space with norm $\Vert\cdot\Vert_{X}$, and
$B:X\times X\rightarrow X$ be a continuous bilinear map, i.e., there exists
$K>0$ such that
\[
\Vert B(x_{1},x_{2})\Vert_{X}\leq K\text{ }\Vert x_{1}\Vert_{X}\text{ }\Vert
x_{2}\Vert_{X},
\]
for all $x_{1},x_{2}\in X$. Given $0<\varepsilon<\frac{1}{4K}$ and $y\in X$
such that $\Vert y\Vert_{X}\leq\varepsilon$, there exists a solution $x\in X$
for the equation $x=y+B(x,x)$ such that $\Vert x\Vert_{X}\leq2\varepsilon$.
The solution $x$ is unique in the closed ball~$\left\{  x\in X:\left\Vert
x\right\Vert _{X}\leq2\varepsilon\right\}  .$ Moreover, the solution depends
continuously on $y$ in the following sense: If $\Vert\tilde{y}\Vert_{X}%
\leq\varepsilon$, $\tilde{x}=\tilde{y}+B(\tilde{x},\tilde{x})$, and
$\Vert\tilde{x}\Vert_{X}\leq2\varepsilon$, then
\[
\Vert x-\tilde{x}\Vert_{X}\leq\frac{1}{1-4K\varepsilon}\Vert y-\tilde{y}%
\Vert_{X}.
\]

\end{lemma}

\begin{remark}
\label{rem-seq}(Picard sequence) The solution given by Lemma \ref{genlem} can
be obtained as the limit in $X$ of the Picard sequence $\{x_{k}\}_{k\in
\mathbb{N}}$ where $x_{1}=y$ and $x_{k+1}=y+B(x_{k},x_{k}),$ for all
$k\in\mathbb{N}$. Moreover, $\left\Vert x_{k}\right\Vert _{X}\leq2\varepsilon$
for all $k\in\mathbb{N}$.
\end{remark}

The following proposition shows that (\ref{dase})-(\ref{ujotas}) is locally in
time well-posed for $L^{\frac{n}{2\gamma-\beta}}$-data.

\begin{proposition}
\label{Proplocal-1}(Local in time solutions) Assume (\ref{cond-1}) and let $q$
be such that
\[
\max\left\{  \frac{\beta-1}{n},\frac{\gamma-1}{n}\right\}  <\frac{1}{q}%
<\min\left\{  \frac{2\gamma-\beta}{n},\frac{n+\beta-2\gamma}{n},\frac
{n+\beta-1}{2n}\right\}  .
\]
If $\theta_{0}\in L^{\frac{n}{2\gamma-\beta}}(\mathbb{R}^{n})$ then there
exists $T>0$ such that (\ref{dase})-(\ref{ujotas}) has a unique mild solution
$\theta$ in the class
\begin{equation}
t^{\eta_{q}}\theta\in BC((0,T);\dot{H}_{q}^{\beta-1}(\mathbb{R}^{n}))\text{
and }\lim_{t\rightarrow0^{+}}t^{\eta_{q}}\left\Vert \theta\right\Vert
_{\dot{H}_{q}^{\beta-1}}=0. \label{loc-sol-1}%
\end{equation}
Moreover, $\theta\in BC([0,T);L^{\frac{n}{2\gamma-\beta}}(\mathbb{R}^{n})).$
\end{proposition}

\textbf{Proof. }\bigskip

\textit{Step 1:} For $T>0$, let us define the Banach space
\[
\mathcal{E}_{T}=\left\{  \theta\text{ measurable};\text{ }t^{\eta_{q}}%
\theta\in BC((0,T);\dot{H}_{q}^{\beta-1}(\mathbb{R}^{n}))\right\}
\]
with norm given by
\begin{equation}
\left\Vert \theta\right\Vert _{\mathcal{E}_{T}}=\sup_{0<t<T}t^{\eta_{q}}%
\Vert\theta(\cdot,t)\Vert_{\dot{H}_{q}^{\beta-1}}. \label{norm1}%
\end{equation}
Due to (\ref{zeroPrincipal}) in Lemma \ref{sol fund} and $\theta_{0}\in
L^{\frac{n}{2\gamma-\beta}}(\mathbb{R}^{n}),$ for any $\varepsilon>0$, there
exists a $T>0$ such that
\begin{equation}
\sup_{0<t<T}t^{\eta_{q}}\Vert G_{\gamma}(t)\theta_{0}\Vert_{\dot{H}_{q}%
^{\beta-1}}\leq\varepsilon. \label{aux-lin}%
\end{equation}
Take $0<\varepsilon<\frac{1}{4K_{4}}$ where $K_{4}$ is as in (\ref{bili-4})
with $r=q$. In view of (\ref{aux-lin}) and (\ref{bili-4}), we can apply Lemma
\ref{genlem} in $\mathcal{E}_{T}$ to obtain a unique solution $\theta(x,t)$
for (\ref{mild}) such that%
\begin{equation}
\ \sup_{0<t<T}t^{\eta_{q}}\Vert\theta\Vert_{\dot{H}_{q}^{\beta-1}}%
\leq2\varepsilon. \label{aux-prop1}%
\end{equation}

Using (\ref{zeroPrincipal}) and an induction argument, one can shows that
every element $\theta_{k}$ of the Picard sequence
\begin{align}
\theta_{1}(x,t)  &  =G_{\gamma}(t)\theta_{0}(x),\label{seq1}\\
\theta_{k+1}(x,t)  &  =\theta_{1}(x,t)+B(\theta_{k},\theta_{k}),\text{ }%
k\in\mathbb{N}, \label{seq2}%
\end{align}
satisfies $\lim_{t\rightarrow0^{+}}t^{\eta_{q}}\left\Vert \theta
_{k}\right\Vert _{\dot{H}_{q}^{\beta-1}}=0.$ Then the second property in
(\ref{loc-sol-1}) follows from the fact that the fixed point $\theta$\ is the
limit in $\mathcal{E}_{T}$ of $\{\theta_{k}\}_{k\in\mathbb{N}}$ (see Remark
\ref{rem-seq}). Further details are left to the reader.

\textit{Step 2: }In what follows we show that $\theta\in$\textbf{
}$BC([0,T);L^{\frac{n}{2\gamma-\beta}}(\mathbb{R}^{n}))$.\ We have that the
recursive sequence (\ref{seq1})-(\ref{seq2}) satisfies (see Remark
\ref{rem-seq})
\begin{equation}
\sup_{0<t<T}t^{\eta_{q}}\Vert\theta_{k}\Vert_{\dot{H}_{q}^{\beta-1}}%
\leq2\varepsilon,\text{ for all }k\in\mathbb{N}. \label{aux-seq-1}%
\end{equation}
Using Lemma \ref{sol fund}, (\ref{bili-2}) with $p=\frac{n}{2\gamma-\beta}$,
and (\ref{aux-seq-1}), we get
\[
\sup_{0<t<T}\Vert\theta_{1}(t)\Vert_{\frac{n}{2\gamma-\beta}}\leq C\Vert
\theta_{0}\Vert_{\frac{n}{2\gamma-\beta}}\text{,}%
\]
and
\begin{align}
\sup_{0<t<T}\Vert\theta_{k+1}(t)\Vert_{\frac{n}{2\gamma-\beta}}  &  \leq
C\Vert\theta_{0}\Vert_{\frac{n}{2\gamma-\beta}}+K_{2}\sup_{0<t<T}t^{\eta_{q}%
}\Vert\theta_{k}(t)\Vert_{\dot{H}_{q}^{\beta-1}}\sup_{0<t<T}\Vert\theta
_{k}(t)\Vert_{\frac{n}{2\gamma-\beta}}\nonumber\\
&  \leq C\Vert\theta_{0}\Vert_{\frac{n}{2\gamma-\beta}}+2\varepsilon K_{2}%
\sup_{0<t<T}\Vert\theta_{k}(t)\Vert_{\frac{n}{2\gamma-\beta}},\text{ for all
}k\in\mathbb{N}. \label{aux-seq-3}%
\end{align}
By reducing $T>0$ in (\ref{aux-lin}) if necessary, we can consider
$0<\varepsilon<\min\{\frac{1}{4K_{4}},\frac{1}{2K_{2}}\}.$ Since
$2K_{2}\varepsilon<1,$ an induction argument in (\ref{aux-seq-3}) shows that
$\{\theta_{k}\}_{k\in\mathbb{N}}$ is uniformly bounded in $L^{\infty
}((0,T);L^{\frac{n}{2\gamma-\beta}}(\mathbb{R}^{n}))$ and then there exists a
subsequence of $\{\theta_{k}\}_{k\in\mathbb{N}}$ (denoted in the same way)
that converges toward $\widetilde{\theta}$ weak$-\ast$ in that space and
consequently in $\mathcal{D}^{\prime}(\mathbb{R}^{n}\times\lbrack0,T)).$
Because $\theta_{k}\rightarrow\theta$ in $\mathcal{E}_{T}$, which implies
convergence in $\mathcal{D}^{\prime}(\mathbb{R}^{n}\times\lbrack0,T)),$ the
uniqueness of the limit in the sense of distributions yields $\theta
=\widetilde{\theta}\in L^{\infty}((0,T);L^{\frac{n}{2\gamma-\beta}}%
(\mathbb{R}^{n}))$. The time-continuity of $\theta$ follows from standard
arguments by using that $\theta$ verifies (\ref{mild}), $\theta\in L^{\infty
}((0,T);L^{\frac{n}{2\gamma-\beta}}(\mathbb{R}^{n}))\cap\mathcal{E}_{T}$, and
the second property in (\ref{loc-sol-1}) (see e.g. \cite{Ka1, Ka2}).

{\vskip 0pt \hfill\hbox{\vrule height 5pt width 5pt depth 0pt} \vskip 12pt}

In the next proposition we investigate the $L^{1}$-persistence of the
solutions obtained in Proposition \ref{Proplocal-1}.

\begin{proposition}
\label{Proplocal-2}Under hypotheses of Proposition \ref{Proplocal-1}. There
exists $T>0$ such that the solution $\theta$ belongs to $BC\left(
[0,T);L^{1}(\mathbb{R}^{n})\right)  $ provided that $\theta_{0}\in
L^{1}(\mathbb{R}^{n})\cap L^{\frac{n}{2\gamma-\beta}}(\mathbb{R}^{n})$.
\end{proposition}

\textbf{Proof.} Let $q$ be such that $1<q^{\prime}<\frac{n}{2\gamma-\beta}.$
From interpolation, we have that $\theta_{0}\in L^{q^{\prime}}(\mathbb{R}%
^{n}).$ Employing (\ref{aux-seq-1}) and the estimate (\ref{bili-2}) with
$p=q^{\prime}$ , we get
\[
\sup_{0<t<T}\Vert\theta_{1}(t)\Vert_{q^{\prime}}\leq C\Vert\theta_{0}%
\Vert_{q^{\prime}}\text{,}%
\]
and
\begin{align*}
\sup_{0<t<T}\Vert\theta_{k+1}(t)\Vert_{q^{\prime}}  &  \leq C\Vert\theta
_{0}\Vert_{q^{\prime}}+K_{2}\sup_{0<t<T}t^{\eta_{q}}\Vert\theta_{k}%
(t)\Vert_{\dot{H}_{q}^{\beta-1}}\sup_{0<t<T}\Vert\theta_{k}(t)\Vert
_{q^{\prime}}\\
&  \leq C\Vert\theta_{0}\Vert_{q^{\prime}}+2K_{2}\varepsilon\sup_{0<t<T}%
\Vert\theta_{k}(t)\Vert_{q^{\prime}},\text{ for all }k\in\mathbb{N}.
\end{align*}
Again reducing $T>0$ if necessary, we can consider $2K_{2}\varepsilon<1$ and
proceed similarly to the end of the proof of Proposition \ref{Proplocal-1} to
obtain that
\begin{equation}
\theta\in BC\left(  [0,T);L^{q^{\prime}}(\mathbb{R}^{n})\right)  .
\label{loc-sol-2}%
\end{equation}
Now we use (\ref{est linear2}), (\ref{bili-3}), (\ref{loc-sol-1}),
(\ref{loc-sol-2}) to estimate
\begin{align*}
\sup_{0<t<T}\Vert\theta(t)\Vert_{1}  &  \leq C\Vert\theta_{0}\Vert_{1}%
+\sup_{0<t<T}\left\Vert B(\theta,\theta)\right\Vert _{1}\\
&  \leq C\Vert\theta_{0}\Vert_{1}+K_{3}T^{1-\frac{1}{2\gamma}-\eta_{q}}%
\sup_{0<t<T}t^{\eta_{q}}\Vert\theta\Vert_{\dot{H}_{q}^{\beta-1}}\,\sup
_{0<t<T}\Vert\theta\Vert_{q^{\prime}}<\infty,
\end{align*}
as required.

{\vskip 0pt \hfill\hbox{\vrule height 5pt width 5pt depth 0pt} \vskip 12pt}

The existence time $T$ in Propositions \ref{Proplocal-1} and \ref{Proplocal-2}
may depend on index $q.$ In the next proposition we show that indeed one can
take a same small time $T>0$ for all $q.$

\begin{proposition}
\label{extnormas} Under hypotheses of Proposition \ref{Proplocal-1}. Let
$\theta$ be the solution given by Proposition \ref{Proplocal-1} with data
$\theta_{0}\in L^{\frac{n}{2\gamma-\beta}}(\mathbb{R}^{n})$. There is a $T>0$
such that
\begin{align}
t^{\eta_{r}}\theta &  \in BC((0,T);\dot{H}_{r}^{\beta-1}(\mathbb{R}%
^{n})),\text{ for all }\frac{n}{2\gamma-\beta}<r<\infty,\label{solutionL}\\
\text{ }t^{\tilde{\eta}_{r}}\theta &  \in BC((0,T);L^{r}(\mathbb{R}%
^{n})),\text{ for all }\frac{n}{2\gamma-\beta}<r<\infty. \label{solutionL2}%
\end{align}

\end{proposition}

\textbf{Proof.} Let $q$ be fixed and as in Proposition \ref{Proplocal-1}.
Given $\frac{n}{2\gamma-\beta}<r<\infty$ verifying $\frac{2\gamma-2}{n}%
-\frac{1}{q}<\frac{1}{r}<\frac{n+\beta-1}{n}-\frac{1}{q},$ we can use
(\ref{bili-4}) instead of (\ref{bili-1}) and proceed just like as in step 2 of
the proof of Proposition \ref{Proplocal-1} to obtain (reducing $T>0$ if
necessary)
\begin{equation}
t^{\eta_{r}}\theta\in BC((0,T);\dot{H}_{r}^{\beta-1}(\mathbb{R}^{n})).
\label{aux-H1}%
\end{equation}
Now let $\frac{2\gamma-2}{n}-\frac{1}{q}<\frac{1}{r}<\frac{2\gamma-1}{n}%
-\frac{1}{q},$ and consider $r<\tilde{r}<\infty$ . Taking $\frac{1}{h}%
=\frac{1}{q}+\frac{1}{z}=\frac{1}{q}+\frac{1}{r}-\frac{\beta-1}{n}$ and
$\delta=\frac{n}{h}-\frac{n}{\tilde{r}}$, it follows that $\delta>0,$
$\frac{\beta+\delta}{2\gamma}<1,$ and $\eta_{q}+\eta_{r}<1$. Then, we can
estimate
\begin{align}
\Vert\theta\Vert_{\dot{H}_{\tilde{r}}^{\beta-1}}  &  \leq Ct^{-\eta_{\tilde
{r}}}\Vert\theta_{0}\Vert_{\frac{n}{2\gamma-\beta}}+\int_{0}^{t}\Vert
G_{\gamma}(t-s)\nabla_{x}\cdot(P[\theta]\theta)(s)\Vert_{\dot{H}_{\tilde{r}%
}^{\beta-1}}\text{ }ds\nonumber\\
&  \leq Ct^{-\eta_{\tilde{r}}}\Vert\theta_{0}\Vert_{\frac{n}{2\gamma-\beta}%
}+C\int_{0}^{t}\Vert G_{\gamma}(t-s)\nabla_{x}\cdot(P[\theta]\theta
)(s)\Vert_{\dot{H}_{h}^{\beta-1+\delta}}\text{ }ds\nonumber\\
&  \leq Ct^{-\eta_{\tilde{r}}}\Vert\theta_{0}\Vert_{\frac{n}{2\gamma-\beta}%
}+C\int_{0}^{t}(t-s)^{-\frac{\beta+\delta}{2\gamma}}\Vert\nabla_{x}%
\cdot(P[\theta]\theta)(s)\Vert_{\dot{H}_{h}^{-1}}\text{ }ds,
\label{aux-proof2}%
\end{align}
where above we have used Sobolev embedding and afterwards (\ref{est linear}).

Now we employ (\ref{fjotas}), H\"{o}lder inequality and Sobolev embedding in
order to estimate%

\begin{align}
&  \text{R.H.S. of (\ref{aux-proof2})}\nonumber\\
&  \leq Ct^{-\eta_{\tilde{r}}}\Vert\theta_{0}\Vert_{\frac{n}{2\gamma-\beta}%
}+C\int_{0}^{t}(t-s)^{-\frac{\beta+\delta}{2\gamma}}\Vert P[\theta
]\theta(s)\Vert_{h}\text{ }ds\nonumber\\
&  \leq Ct^{-\eta_{\tilde{r}}}\Vert\theta_{0}\Vert_{\frac{n}{2\gamma-\beta}%
}+C\int_{0}^{t}(t-s)^{-\frac{\beta+\delta}{2\gamma}}\Vert P[\theta]\Vert
_{q}\Vert\theta(s)\Vert_{z}\text{ }ds\nonumber\\
&  \leq Ct^{-\eta_{\tilde{r}}}\Vert\theta_{0}\Vert_{\frac{n}{2\gamma-\beta}%
}+C\int_{0}^{t}(t-s)^{-\frac{\beta+\delta}{2\gamma}}\Vert\theta\Vert_{\dot
{H}_{q}^{\beta-1}}\Vert\theta\Vert_{\dot{H}_{r}^{\beta-1}}\text{
}ds\nonumber\\
&  \leq Ct^{-\eta_{\tilde{r}}}\Vert\theta_{0}\Vert_{\frac{n}{2\gamma-\beta}%
}+C\int_{0}^{t}(t-s)^{-\frac{\beta+\delta}{2\gamma}}s^{-\eta_{q}-\eta_{r}%
}\text{ }ds\left(  \sup_{0<t<T}t^{\eta_{q}}\Vert\theta\Vert_{\dot{H}%
_{q}^{\beta-1}}\right)  \left(  \sup_{0<t<T}t^{\eta_{r}}\Vert\theta\Vert
_{\dot{H}_{r}^{\beta-1}}\right) \label{aux-T}\\
&  \leq Ct^{-\eta_{\tilde{r}}}\Vert\theta_{0}\Vert_{\frac{n}{2\gamma-\beta}%
}+Ct^{-\eta_{\tilde{r}}}\int_{0}^{1}(1-s)^{-\frac{\beta+\delta}{2\gamma}%
}s^{-\eta_{q}-\eta_{r}}\text{ }ds\nonumber\\
&  \leq Ct^{-\eta_{\tilde{r}}},\nonumber
\end{align}
and then we arrive at (\ref{aux-H1}) with $\tilde{r}$ in place of $r$, and
with the same existence time $T>0$. From interpolation, notice that
(\ref{aux-H1}) also holds true for every $r=l$ such that $\frac{n}%
{2\gamma-\beta}<l<\tilde{r}.$ Since $\tilde{r}>r$ is arbitrary, we obtain
(\ref{aux-H1}) with $r=l$ ( and the same $T>0$), for all $\frac{n}%
{2\gamma-\beta}<l<\infty,$ which gives (\ref{solutionL}).

The proof of (\ref{solutionL2}) can be performed in a similar way by using
(\ref{bili-1}) instead of (\ref{bili-4}). {\vskip 0pt \hfill
\hbox{\vrule height 5pt width 5pt depth 0pt} \vskip 12pt}

\bigskip

\subsection{Proof of Theorem \ref{teoglobal}}

\bigskip

\subsubsection{Step 1: \ Local smoothness and maximum principle}

\hspace{0.55cm}The solutions obtained in Proposition \ref{Proplocal-1} are
instantaneously $C^{\infty}$-smoothed for any $t>0$ belonging to the existence
interval $(0,T).\,$\ This smooth effect holds for several parabolic equations
in several frameworks, like e.g. $L^{p},$ weak-$L^{p}$, Morrey, Besov-Morrey,
when mild solutions are constructed by using time-weighted norms of Kato type
(see \cite{Ka1}). Precisely, adapting arguments from \cite{Ka1} (see also
\cite{Car-Fer2}), one can obtain that the solution verifies%
\begin{equation}
\partial_{t}^{m}\nabla_{x}^{k}\theta(x,t)\in C\left(  (0,T);L^{\frac
{n}{2\gamma-\beta}}(\mathbb{R}^{n})\cap L^{q}(\mathbb{R}^{n})\right)  ,
\label{derivadas}%
\end{equation}
for all $\frac{n}{2\gamma-\beta}<q<\infty,$ $m\in\{0\}\cup\mathbb{N}$ and
multi-index $k\in(\{0\}\cup\mathbb{N)}^{n},$ where $T>0$ is the existence time
given in Proposition \ref{extnormas}. In particular, it follows that
$\theta(x,t)\in C^{\infty}(\mathbb{R}^{n}\times(0,T))$ and $\theta(t)\in
L^{\infty}(\mathbb{R}^{n})$ with
\begin{equation}
\left\Vert \theta(t)\right\Vert _{\infty}\leq C\left\Vert \theta(t)\right\Vert
_{\frac{n^{2}}{2\gamma-\beta}}^{\alpha}\left\Vert \nabla_{x}\theta
(t)\right\Vert _{\frac{n^{2}}{2\gamma-\beta}}^{1-\alpha},\text{ }
\label{Aux-inf}%
\end{equation}
for all $0<t<T,$ where $\alpha=\frac{n+\beta-2\gamma}{n}$. If further
$\theta_{0}\in L^{1}(\mathbb{R}^{n})\cap L^{\frac{n}{2\gamma-\beta}%
}(\mathbb{R}^{n})$ then $q$ in (\ref{derivadas}) can be taken in the range
$1<q<\infty.$

Due to (\ref{derivadas}) we have that $\theta$ verifies (\ref{dase}%
)-(\ref{ujotas}) in the classical sense and $\partial_{t}^{m}\nabla_{x}%
^{k}\theta(x,t)\rightarrow0$ when $|x|\rightarrow\infty,$ for all $0<t<T.$ In
view of $\nabla\cdot$ $P[\theta]=0$, we can integrate by parts to obtain
\begin{align}
\frac{\partial}{\partial t}\Vert\theta(t)\Vert_{p}^{p}  &  =p\int
_{\mathbb{R}^{n}}\theta(t)^{p-1}\frac{\partial}{\partial t}\theta
(t)dx\nonumber\\
&  =p\int_{\mathbb{R}^{n}}\theta(t)^{p-1}\left(  -(-\Delta)^{\gamma}%
\theta-\nabla_{x}\cdot(P[\theta]\theta)\right)  dx\nonumber\\
&  =-p\int_{\mathbb{R}^{n}}\theta(t)^{p-1}(-\Delta)^{\gamma}\theta dx\leq
-\int_{\mathbb{R}^{2}}\left\vert (-\Delta)^{\frac{\gamma}{2}}(\theta^{\frac
{p}{2}})\right\vert ^{2}dx, \label{ener}%
\end{align}
for all $t\in(0,T)$. The last inequality in (\ref{ener}) can be found in
\cite{Const3, Cordoba1} (see also \cite{Ju}).

In view of the estimate (\ref{ener}), we have that $L^{p}$-norms of
$\theta(t)$ are non-increasing in $(0,T).$ If $\theta_{0}\in L^{\frac
{n}{2\gamma-\beta}}(\mathbb{R}^{n})$ and $\theta_{0}\in L^{1}(\mathbb{R}%
^{n})\cap L^{\frac{n}{2\gamma-\beta}}(\mathbb{R}^{n}),$ we obtain
respectively
\begin{equation}
\Vert\theta(t)\Vert_{\frac{n}{2\gamma-\beta}}\leq\Vert\theta(t_{0}%
)\Vert_{\frac{n}{2\gamma-\beta}}\text{ and }\Vert\theta(t)\Vert_{1}\leq
\Vert\theta(t_{0})\Vert_{1}, \label{aux-max-1}%
\end{equation}
for $0<t_{0}\leq t<T.$

Making $t_{0}\rightarrow0^{+}$ in (\ref{aux-max-1}), it follows that the
solution $\theta(x,t)$ satisfies
\begin{equation}
\Vert\theta(t)\Vert_{\frac{n}{2\gamma-\beta}}\leq\Vert\theta_{0}\Vert
_{\frac{n}{2\gamma-\beta}}\text{ and }\Vert\theta(t)\Vert_{1}\leq\Vert
\theta_{0}\Vert_{1}, \label{max-prin}%
\end{equation}
for all $t\in(0,T),$ when $\theta_{0}\in L^{\frac{n}{2\gamma-\beta}%
}(\mathbb{R}^{n})$ and $\theta_{0}\in L^{1}(\mathbb{R}^{n})\cap L^{\frac
{n}{2\gamma-\beta}}(\mathbb{R}^{n}),$ respectively.

\subsubsection{\bigskip Step 2: Extension of the local solution}

\hspace{0.55cm}We start by making the following observation: if $\frac
{n}{2\gamma-\beta}<q<\infty$ and $\theta_{0}\in L^{q}(\mathbb{R}^{n})$ then
\begin{equation}
t^{\eta_{q}}\Vert G_{\gamma}(t)\theta_{0}\Vert_{\dot{H}_{q}^{\beta-1}}\leq
Ct^{\tilde{\eta}_{q}}\Vert\theta_{0}\Vert_{q}\rightarrow0\text{ when
}t\rightarrow0^{+}. \label{aux-ext-1}%
\end{equation}
Therefore, for $\theta_{0}\in L^{\frac{n}{2\gamma-\beta}}(\mathbb{R}^{n})\cap
L^{q}(\mathbb{R}^{n})$, the existence time $T>0$ obtained in Proposition
\ref{Proplocal-1} can be taken depending on the norm $\Vert\theta_{0}\Vert
_{q}.$ Indeed it can be chosen as
\begin{equation}
T=\left(  \frac{\varepsilon}{C\Vert\theta_{0}\Vert_{q}}\right)  ^{\frac
{1}{\tilde{\eta}_{q}}}, \label{aux-T-1}%
\end{equation}
where $0<\varepsilon<\frac{1}{4K_{4}}$ and $C>0$ is as in (\ref{aux-ext-1}).

Now let $\theta_{0}\in L^{\frac{n}{2\gamma-\beta}}(\mathbb{R}^{n})$. From
Proposition \ref{Proplocal-1}, there exist $M_{0}>0$, $T_{0}>0$ and a unique
mild solution for (\ref{dase})-(\ref{ujotas}) in $(0,T_{0})$ such that
\begin{equation}
\sup_{0<t<T_{0}}t^{\eta_{q}}\Vert\theta(t)\Vert_{\dot{H}_{q}^{\beta-1}}\leq
M_{0}\text{ and }\sup_{0<t<T_{0}}\Vert\theta(t)\Vert_{\frac{n}{2\gamma-\beta}%
}\leq\Vert\theta_{0}\Vert_{\frac{n}{2\gamma-\beta}}, \label{aux-proof1}%
\end{equation}
where (\ref{max-prin}) has been used in the second inequality in
(\ref{aux-proof1}).

Let us now denote
\begin{equation}
T=\sup\left\{  \tilde{T}>0;\theta\in C((0,\tilde{T});\dot{H}_{q}^{\beta
-1}(\mathbb{R}^{n})\cap L^{\frac{n}{2\gamma-\beta}}(\mathbb{R}^{n}%
)),\sup_{0<t<\tilde{T}}t^{\eta_{q}}\Vert\theta\Vert_{\dot{H}_{q}^{\beta-1}%
}<\infty,\sup_{0<t<\tilde{T}}\Vert\theta\Vert_{\frac{n}{2\gamma-\beta}}%
\leq\Vert\theta_{0}\Vert_{\frac{n}{2\gamma-\beta}}\right\}  . \label{aux-def1}%
\end{equation}
We desire to prove that $T=\infty.$ Suppose by contradiction that $T<\infty$,
and let $a=\theta(T-\varepsilon)$ where $0<\varepsilon<T$ will be choose
later. Proposition \ref{extnormas} gives that $a\in L^{\frac{n}{2\gamma-\beta
}}(\mathbb{R}^{n})\cap L^{q}(\mathbb{R}^{n}),$ for all $\frac{n}{2\gamma
-\beta}<q<\infty$. Moreover, if $\varepsilon<\frac{T}{2}$ we get $\Vert
\theta(T-\varepsilon)\Vert_{q}\leq\Vert\theta(\frac{T}{2})\Vert_{q}$.
Therefore, taking $\varepsilon<T/2$ and $a$ as initial data, we have that
given $0<M_{1}<\frac{1}{4K_{4}}$ there exist $T_{1}>0$ and a unique mild
solution $\tilde{\theta}$ for (\ref{dase})-(\ref{ujotas}) such that%
\[
\tilde{\theta}\in C((T-\varepsilon,T_{1}+T-\varepsilon);L^{\frac{n}%
{2\gamma-\beta}}(\mathbb{R}^{n})\cap\dot{H}_{q}^{\beta-1}(\mathbb{R}^{n}))
\]
and
\begin{align}
\sup_{T-\varepsilon<t<T_{1}+T-\varepsilon}(t-(T-\varepsilon))^{\eta_{q}%
}\left\Vert \tilde{\theta}(t)\right\Vert _{\dot{H}_{q}^{\beta-1}}  &
\leq2M_{1}\label{aux-ext-10}\\
\sup_{T-\varepsilon<t<T_{1}+T-\varepsilon}\Vert\tilde{\theta}(t)\Vert
_{\frac{n}{2\gamma-\beta}}  &  \leq\Vert\theta(T-\varepsilon)\Vert_{\frac
{n}{2\gamma-\beta}}.\nonumber
\end{align}
From uniqueness part of Proposition \ref{Proplocal-1}, it follows that
$\theta=\tilde{\theta}$ in $(T-\varepsilon,T).$ In view of (\ref{aux-T-1}), we
can choose $T_{1}=\min\left\{  \left(  \frac{M_{1}}{C\Vert\theta(\frac{T}%
{2})\Vert_{q}}\right)  ^{\frac{1}{\tilde{\eta}_{q}}},T\right\}  $. Taking
$0<\varepsilon<\min\{\frac{T}{2},T_{1}\}$ and $T_{2}=T_{1}+T-\varepsilon,$ we
have that $T<T_{2}$ and get a solution%
\[
\theta\in C((0,T_{2});L^{\frac{n}{2\gamma-\beta}}(\mathbb{R}^{n})\cap\dot
{H}_{q}^{\beta-1}(\mathbb{R}^{n}))
\]
such that
\[
\sup_{0<t<\tilde{T}}t^{\eta_{q}}\left\Vert \theta(t)\right\Vert _{\dot{H}%
_{q}^{\beta-1}}<\infty\text{ and }\sup_{0<t<\tilde{T}}\Vert\theta
(t)\Vert_{\frac{n}{2\gamma-\beta}}\leq\Vert\theta_{0}\Vert_{\frac{n}%
{2\gamma-\beta}},
\]
for all $0<\tilde{T}<T_{2},$ which contradicts the maximality of $T$ in
(\ref{aux-def1}). Therefore $T=\infty$ and we are done.

\subsubsection{Step 3: Global $L^{q}$-decay of solutions}

\bigskip\hspace{0.55cm}We will prove only the part of the statement
corresponding to the case $\theta_{0}\in L^{\frac{n}{2\gamma-\beta}%
}(\mathbb{R}^{n})$. The estimate (\ref{Decay-Lq-2}) for $\theta_{0}\in
L^{1}(\mathbb{R}^{n})\cap L^{\frac{n}{2\gamma-\beta}}(\mathbb{R}^{n})$ follows
similarly to the first one by using $\left\Vert \theta(t)\right\Vert _{1}%
\leq\left\Vert \theta_{0}\right\Vert _{1}$ and the sequence $q_{k}=2^{k}$
instead of $\left\Vert \theta(t)\right\Vert _{\frac{n}{2\gamma-\beta}}%
\leq\left\Vert \theta_{0}\right\Vert _{\frac{n}{2\gamma-\beta}}$ and
$q_{k}=\frac{n}{2\gamma-\beta}2^{k}$.

Since we have extended the solution $\theta$, it follows that (\ref{derivadas}%
) and (\ref{Aux-inf}) hold true for $T=\infty$. Then
\begin{equation}
\left\Vert \theta(\cdot,t)\right\Vert _{\infty}<\infty,\text{ for all }t>0.
\label{Aux-inf-2}%
\end{equation}

Now we proceed as in \cite{Car-Fer1} and \cite{Ka2}. In view of the
Gagliardo-Nirenberg inequality, we have that
\begin{equation}
\left\Vert \phi\right\Vert _{2}\leq C\left\Vert \phi\right\Vert _{1}^{\alpha
}\left\Vert (-\Delta)^{\frac{\gamma}{2}}\phi\right\Vert _{2}^{1-\alpha}\text{
with }\alpha=\frac{2\gamma}{n+2\gamma}. \label{Ga-Ni}%
\end{equation}
Taking $\phi=\theta^{\frac{q}{2}}$ in (\ref{Ga-Ni}), it follows that \
\begin{equation}
\left\Vert \theta\right\Vert _{q}^{q(\frac{n+2\gamma}{n})}\leq C\left\Vert
\theta\right\Vert _{\frac{q}{2}}^{\frac{2q\gamma}{n}}\left\Vert (-\Delta
)^{\frac{\gamma}{2}}(\theta^{\frac{q}{2}})\right\Vert _{2}^{2}\text{.}
\label{aux-Global1}%
\end{equation}
Denoting $\psi_{q}(t)=$ $\Vert\theta(t)\Vert_{q}^{q}$ , we obtain from
(\ref{ener}) and (\ref{aux-Global1}) that
\begin{equation}
\frac{\partial}{\partial t}\psi_{q}\leq-C(\psi_{\frac{q}{2}})^{-\frac{4\gamma
}{n}}\psi_{q}^{\frac{n+2\gamma}{n}}. \label{ineq}%
\end{equation}
The differential inequality (\ref{ineq}) can be solved by an induction
procedure. In fact, using the first inequality in (\ref{max-prin}) and
considering the sequence $q_{k}=\frac{n}{2\gamma-\beta}2^{k}$ for $k\geq0,$ we
arrive at
\begin{equation}
\psi_{q_{k}}(t)\leq M_{q_{k}}t^{-\frac{n}{2}(\frac{2^{k}-1}{\gamma})}\,,\,
\label{mx}%
\end{equation}
where
\[
M_{q_{0}}=\left\Vert \theta_{0}\right\Vert _{\frac{n}{2\gamma-\beta}}\text{
and }\,M_{q_{k}}=\left(  \frac{n(2^{k}-1)}{2C\gamma}\right)  ^{\frac
{n}{2\gamma}}M_{\frac{q_{k}}{2}}^{2}\text{, \ \ for }k\in\mathbb{N}.
\]
It follows that
\begin{align*}
M_{q_{k}}^{\frac{1}{q_{k}}}  &  =\left(  \frac{n(2^{k}-1)}{2C\gamma}\right)
^{\frac{n}{2\gamma q_{k}}}M_{q_{k-1}}^{\frac{1}{q_{k-1}}}=\left(
\frac{n(2^{k}-1)}{2C\gamma}\right)  ^{\frac{n}{2\gamma2^{k}q_{0}}}\left(
\frac{n(2^{k-1}-1)}{2C\gamma}\right)  ^{\frac{n}{2\gamma2^{k-1}q_{0}}%
}M_{q_{k-2}}^{\frac{1}{q_{k-2}}}=...\\
&  =\left[  \Pi_{i=1}^{k}\left(  \frac{n(2^{i}-1)}{2C\gamma}\right)
^{\frac{n}{2^{i}2\gamma q_{0}}}\right]  (M_{q_{0}})^{\frac{1}{q_{0}}},\text{
for all }k\in\mathbb{N},
\end{align*}
and then
\begin{align}
\Vert\theta(t)\Vert_{2^{k}q_{0}}  &  \leq M_{q_{k}}^{\frac{1}{q_{k}}%
}\,t^{-\frac{n}{2}(\frac{2^{k}-1}{\gamma q_{k}})}\nonumber\\
&  =\left(  \Pi_{i=1}^{k}\left(  \frac{n}{2}\frac{2^{i}-1}{C\gamma}\right)
^{\frac{n}{2^{i}2\gamma q_{0}}}\right)  (M_{q_{0}})^{\frac{1}{q_{0}}%
}\,t^{-\frac{n}{2\gamma q_{0}}(\frac{2^{k}-1}{2^{k}})}, \label{aux-global3}%
\end{align}
where $q_{0}=\frac{n}{2\gamma-\beta}.$ In view of (\ref{Aux-inf-2}), we can
make $k\rightarrow\infty$ in (\ref{aux-global3}) to obtain
\begin{equation}
\Vert\theta(t)\Vert_{\infty}\leq C\Vert\theta_{0}\Vert_{\frac{n}{2\gamma
-\beta}}^{\frac{2\gamma-\beta}{n}}\,t^{-\frac{2\gamma-\beta}{2\gamma}}.
\label{aux-global4}%
\end{equation}
Interpolating the first inequality in (\ref{max-prin}) with (\ref{aux-global4}%
), the result is%

\[
\Vert\theta(t)\Vert_{q}\leq C\,t^{-\tilde{\eta}_{q}},\text{ for all }\frac
{n}{2\gamma-\beta}\leq q\leq\infty,
\]
as required.

The uniqueness statement follows from the local uniqueness property in
Proposition \ref{Proplocal-1}.

{\vskip 0pt \hfill\hbox{\vrule height 5pt width 5pt depth 0pt} \vskip 12pt}

\subsection{Proof of Theorem \ref{Teo-sym}}

\hspace{0.55cm}\textbf{Part (i): }We will prove only the odd part of the
statement since the even one follows similarly. Let $\theta$ be the solution
of Proposition \ref{Proplocal-1} with existence time $T>0$. From step 2 of the
proof of Theorem \ref{teoglobal}, $\theta$ can be extended by using
Proposition \ref{Proplocal-1} and solving (\ref{dase})-(\ref{ujotas})
consecutively with initial data $\theta(\frac{T}{2}),$ $\theta(\frac{T}%
{2}+T_{1}),$ $\theta(\frac{T}{2}+2T_{1})$ and so on, where $T_{1}%
=\min\{\left(  \frac{\varepsilon}{C\Vert\theta(\frac{T}{2})\Vert_{q}}\right)
^{\frac{1}{\tilde{\eta}_{q}}},T\}$, $\varepsilon=\frac{1}{8K_{4}}$, and $C$ as
in (\ref{aux-ext-1}). Because of that, it is sufficient to show the following
claim: if $\theta_{0}\in L^{\frac{n}{2\gamma-\beta}}$ is odd then so is the
solution $\theta(x,t)$ given by Proposition \ref{Proplocal-1}, for all
$t\in(0,T)$. In fact, notice that we can use this claim repeatedly to show
that the global solution $\theta(x,t)$ is odd, for all $t>0.$

Let $\psi(x,t)=G_{\gamma}(t)\theta_{0}.$\ We have that $\theta_{0}%
(-x)=-\theta_{0}(x)$ is equivalent to
\begin{equation}
-\widehat{\theta_{0}}(\xi)=[\theta_{0}(-x)]^{\wedge}(\xi)=\widehat{\theta_{0}%
}(-\xi)\text{ in }\mathcal{S}^{\prime}(\mathbb{R}^{n})\text{.}
\label{aux-sym1}%
\end{equation}
It follows from (\ref{aux-sym1}) that
\begin{align*}
\lbrack\psi(-x,t)]^{\wedge}(\xi)  &  =e^{-t\left\vert \xi\right\vert
^{2\gamma}}\widehat{\theta_{0}}(-\xi)\\
&  =-e^{-t\left\vert \xi\right\vert ^{2\gamma}}\widehat{\theta_{0}}%
(\xi)=-\widehat{\psi(x,t)}(\xi),
\end{align*}
which shows that $G_{\gamma}(t)\theta_{0}$ is odd, for each fixed $t>0.$ Also,
if $\theta$ is odd then $\nabla\theta$ is even, because
\[
\nabla(\theta(x,t))=\nabla(-\theta(-x,t))=(\nabla\theta)(-x,t)).
\]
Recall that $A=(a_{ij})$ and
\[
P(\xi)=(\widetilde{P}_{1}(\xi),...,\widetilde{P}_{n}(\xi))
\]
where
\[
\widetilde{P}_{j}(\xi)=\sum_{i=1}^{n}a_{ij}\frac{\xi_{i}I}{\left\vert
\xi\right\vert ^{2}}P_{i}(\xi).
\]
It follows that%
\[
\widehat{u}(-\xi)=(\widehat{u_{1}}(-\xi),...,\widehat{u_{n}}(-\xi
))=P(-\xi)\widehat{\theta}(-\xi,t)\text{ \ \ }%
\]
with
\begin{align*}
P(-\xi)  &  =\frac{I}{\left\vert \xi\right\vert ^{2}}[-\xi_{1}P_{1}%
(-\xi),...,-\xi_{n}P_{n}(-\xi)]\text{ }A\\
&  =\frac{I}{\left\vert \xi\right\vert ^{2}}[\xi_{1}P_{1}(\xi),...,\xi
_{n}P_{n}(\xi)]\text{ }A=P(\xi),
\end{align*}
because $P_{i}$'s are odd. Therefore $u=P[\theta]$ is odd when $\theta$ is
odd, and then $(u\cdot\nabla\theta)=(P[\theta]\cdot\nabla\theta)$ is odd too.
Hence if $\theta$ is odd then so is $B(\theta,\theta)$.

So, employing an induction argument, one can prove that each element
$\theta_{k}$ of the Picard sequence (\ref{seq1})-(\ref{seq2}) is odd. Since
$\theta_{k}\rightarrow\theta$ in the norm (\ref{norm1}), then the sequence
(\ref{seq1})-(\ref{seq2}) also converges (up to a subsequence) to $\theta$
a.e. $x\in\mathbb{R}^{n}$, for all $t\in(0,T)$. It follows that $\theta(x,t)$
is odd, for each fixed $t\in(0,T),$ because pointwise convergence preserves
odd symmetry. This shows the desired claim.

\textbf{Part (ii): }From the same reasons given in \textit{part (i)}, we need
only to prove that the local solution of Proposition \ref{Proplocal-1} is
radially symmetric whenever $\theta_{0}$ and $div_{\xi}(P(\xi))$ are too. For
that matter, we first observe that $G_{\gamma}(t)\theta_{0}$ is radial because
$\theta_{0}$ and the kernel $\hat{g}_{\gamma}(\xi,t)=e^{-|\xi|^{2\gamma}t}$
are radial, for all $t>0.$ Also, for $\theta$ radially symmetric, we have
that
\begin{align}
(u\cdot\nabla\theta)  &  =\sum_{j=1}^{n}u_{j}\partial_{x_{j}}\theta
=\frac{\theta^{\prime}(r)}{r}\sum_{j=1}^{n}u_{j}x_{j}\nonumber\\
&  =I\frac{\theta^{\prime}(r)}{r}\sum_{j=1}^{n}(\partial_{\xi_{j}}%
\widehat{u_{j}})^{\vee}=I\frac{\theta^{\prime}(r)}{r}\sum_{j=1}^{n}\left(
\partial_{\xi_{j}}\widetilde{P}_{j}(\xi)\widehat{\theta}\right)  ^{\vee
}\nonumber\\
&  =I\frac{\theta^{\prime}(r)}{r}\left(  \widehat{\theta}(\xi,t)\left(
div_{\xi}(P(\xi)\right)  \right)  ^{^{\vee}}. \label{aux-rad1}%
\end{align}
It follows from (\ref{aux-rad1}) that if $\theta$ and $\left(  div_{\xi}%
(P(\xi)\right)  $ are radial then so is $(u\cdot\nabla\theta)$. Using that
$G_{\gamma}(t)$ preserves radiality, we obtain that $B(\theta,\theta)$ defined
in (\ref{termo bilinear}) is radially symmetric, for each $t\in(0,T)$,
whenever $\theta$ is too. Analogously to \textit{part (i)}, we now can use
induction in order to show that each function $\theta_{k}$ defined in
(\ref{seq1})-(\ref{seq2}) is also radially symmetric. Since $\theta_{k}$
converges (up to a subsequence) to $\theta$ a.e. $x\in$ $\mathbb{R}^{n}$, for
each $t\in(0,T)$, we obtain the required conclusion.

{\vskip 0pt \hfill\hbox{\vrule height 5pt width 5pt depth 0pt}

\end{document}